\documentclass[12pt, reqno]{amsart}
\usepackage{amsmath, amsthm, amscd, amsfonts, amssymb, graphicx, color}
\usepackage{setspace}
\usepackage{mathrsfs}
\usepackage{multicol}
\usepackage{tikz-cd}
\usepackage[bookmarksnumbered, colorlinks, plainpages]{hyperref}
\hypersetup{colorlinks=true,linkcolor=red, anchorcolor=green, citecolor=cyan, urlcolor=red, filecolor=magenta, pdftoolbar=true}

\newtheorem{theorem}{Theorem}[section]
\newtheorem{lemma}[theorem]{Lemma}

\newtheorem{corollary}[theorem]{Corollary}
\theoremstyle{definition}

\theoremstyle{remark}

\numberwithin{equation}{section}

\newcommand{\NN}{\mathbb{N}}

\newcommand{\CC}{\mathbb {C}}

\newcommand{\R}{\mathbb{R}}
\newcommand{\Z}{\mathbb{Z}}

\begin{document}
\setcounter{page}{1}
\title[Weighted superposition operators on Fock spaces  ]{Weighted superposition operators on Fock spaces }
\author[Tesfa  Mengestie]{Tesfa  Mengestie }
\address{ Mathematics Section \\
Western Norway University of Applied Sciences\\
Klingenbergvegen 8, N-5414 Stord, Norway}
\email{Tesfa.Mengestie@hvl.no}
 \begin{abstract}
 We characterize
 all pairs of entire functions $(u,\psi)$  for which the induced weighted superposition operator $S_{(u,\psi)}$ transforms one Fock space into another Fock space.
   Further analytical structures  like boundedness and Lipschitz continuity of $S_{(u,\psi)}$  are described. We, in particular,  show   the Fock spaces
    support no compact weighted  superposition operator.
    \end{abstract}
\maketitle
\section{Introduction}
The theory of superposition  operator has  a long history in the context of real valued functions \cite{JAPP}.  In contrast, there have been only  some studies  on spaces of  analytic functions which includes  Hardy spaces \cite{GAC}, Bergman spaces \cite{GAJ}, Dirichlet type spaces \cite{BFB},  Bloch type  spaces \cite{VV,RC},  and some weighted Banach spaces over the disc \cite{JB}. The goal of this note is to study the operator on the Fock spaces     $\mathcal{F}_p$.  We recall that   $\mathcal{F}_p$  is  the   space of   entire functions $f$ for which
 \begin{align*}
 \|f\|_{ p}= \begin{cases}\Big(\frac{p}{2\pi} \int_{\CC} |f(z)|^p
e^{-\frac{p|z|^2}{2}}  dA(z)\Big)^{\frac{1}{p}} <\infty, \ \ \ \ 0<p<\infty\\
\sup_{z\in \CC} |f(z)|e^{-\frac{|z|^2}{2}}  <\infty, \ \ \ \ p= \infty,
\end{cases}
\end{align*}  where $dA$ is the usual  Lebesgue area measure on the complex plane  $\CC$.
The space $\mathcal{F}_2$ is a reproducing kernel Hilbert space with kernel function $K_w(z)= e^{\overline{w}z}$. For each $w\in \CC$, a calculation shows  the function   $k_w=\|K_w\|_2^{-1}K_w \in \mathcal{F}_p$  and     $\|k_w\|_p= 1 $ for all $p$.  By \cite[p. 37]{Zhu}, for each  entire function $f$ and  $p\neq \infty$
 \begin{align}
\label{global}
 |f(z)|
 \leq e^{\frac{|z|^2}{2}}\bigg( \int_{D(z,1 )} |f(w)|^p
e^{-\frac{p|w|^2}{2}} dA(w)\bigg)^{1/p}\leq \bigg(\frac{2\pi}{p}\bigg)^{\frac{1}{p}}e^{\frac{|z|^2}{2}}  \|f\|_p,
 \end{align} where  $D(z,1 )$ is the disc with center $z$ and radius  $1$. The  same conclusion  with the factor $\big(2\pi/p\big)^{1/p}$ replaced by $1$ holds whenever $p= \infty$.

In this note, we   study  the
 weighted superposition operator on Fock spaces  where the operator $S_\psi$ is covered as a particular case.   Let $(u,\psi)$ be a pair of  holomorphic functions on $\CC$. For two  metric spaces  $\mathcal{ H}_1$ and $\mathcal{ H}_2$, the weighted superposition operator $S_{(u,\psi)}: \mathcal{ H}_1 \to \mathcal{ H}_2$ is defined by   $S_{(u,\psi)}f=  M_u S_\psi(f)$   where $S_\psi f= \psi(f)$ and $M_uf= uf$ are respectively the superposition and  multiplication operators.  Recently,  $S_{(u,\psi)}$ was  studied on the  Bergman and Bloch spaces \cite{DOM}.

The first main question now is to identify  which pairs of analytic symbols $(u, \psi)$  define  weighted superposition operators
$S_{(u,\psi)}$ from  $\mathcal{F}_p$  into  $\mathcal{F}_q$. The questions for $S_{(u,\psi)}$ are in general technically more difficult than the corresponding questions
for $S_{\psi}$ since the presence of the multiplier $u$ can complicate proofs and arguments.   We may first begin with a simple example that illustrates  the problem at hand.   Let $f= \alpha$ be a constant  and $g(z)=z$. Then,   for $0<q\leq \infty$ and  a non-zero $\psi$,
\begin{align*}
\|S_{(u,\psi)}f\|_q= \| \psi(\alpha) u\|_q= |\psi(\alpha) | \| u\|_q \ \ \text{and}\ \ \|S_{(u,\psi)}g\|_q= \| u\psi\|_q.
\end{align*} This shows that if $ S_{(u,\psi)}$ maps $\mathcal{F}_p$  into  $\mathcal{F}_q$, then  both   $u$ and $u\psi$ belong to $\mathcal{F}_q$.  On the other hand, set for example  $q=2$ and consider the function
$h(z)= \sin(\frac{z^2}{2})/z^2.$
 Then  $h$ is an entire function which belongs to  $\mathcal{F}_2$. To see this, observe that
when $|z|= r$ gets larger, then $|h(z)|^2\simeq e^{r^2}/r^4 $ and
\begin{align*}
\int_0^{2\pi}\int_0^\infty |h(re^{it})|^2 r e^{-r^2} dt dr <\infty.
\end{align*}
However, the function $zh$  is not in $\mathcal{F}_2$ since   $|zh(z)|^2r e^{-r^2}\simeq 1/r$
for larger $|z|$. It follows that   $S_{(z,z)}$ fails to map $ \mathcal{F}_2$ into itself   while it is easy to see that  $S_z$ does.  This exhibits the existence of some  degree of interplay between $u$ and $\psi$ in defining $ S_{(u,\psi)}$  on  Fock spaces. Our next main result provides their precise  interplay.
\begin{theorem}\label{thm0}
Let $\psi$ and $u$ be  nonzero entire functions on  $\CC$, and $0<p,q\leq  \infty$.
\begin{enumerate}
\item If $p\leq q$, then  the following statements are equivalent.
\begin{enumerate}
\item $ S_{(u,\psi)}$ maps $\mathcal{F}_p$ into $  \mathcal{F}_q$;
\item  Either  $\psi(z)= az+b$ for some $a, b \in \CC$  and $u$ is a constant or $\psi$ is a constant and $ u\in  \mathcal{F}_q$. If $u$ is in addition  non-vanishing, then
\begin{align}
\label{form101}
 u(z)=u(0) e^{a_1z+a_2z^2} , \ \  a_1, a_2 \in \CC \ \text{ and}\   |a_2|< 1/2;
 \end{align}
\item  $ S_{(u,\psi)}:  \mathcal{F}_p\rightarrow  \mathcal{F}_q$   is bounded;
\item $ S_{(u,\psi)}:  \mathcal{F}_p\rightarrow  \mathcal{F}_q$ is globally Lipschitz continuous.
\end{enumerate}
\item If $p>q$, then the following statements are equivalent.
\begin{enumerate}
\item $ S_{(u,\psi)}$ maps $\mathcal{F}_p$ into $  \mathcal{F}_q$;
\item $\psi$ is a   constant and   $ u\in  \mathcal{F}_q$.  If $u$ is in addition  non-vanishing, then  the representation in  \eqref{form101} holds;
\item  $ S_{(u,\psi)}:  \mathcal{F}_p\rightarrow  \mathcal{F}_q$ is bounded;
\item $ S_{(u,\psi)}:  \mathcal{F}_p\rightarrow  \mathcal{F}_q$ is globally Lipschitz continuous.
\end{enumerate}
\item The map $ S_{(u,\psi)}:  \mathcal{F}_p\rightarrow  \mathcal{F}_q$ cannot be compact for any pair of  $p$ and $q$.
\end{enumerate}
\end{theorem}
The proof of Theorem~\ref{thm0} follows once we prove Lemma~\ref{order},  Theorem~\ref{thm3}, and   Theorem~\ref{thm4}   in the next section.

The operator $S_{(u,\psi)}$ reduces to  $S_{\psi}$ and $M_u$  when   $u=1$ and $\psi(z)= z$ respectively. Consequently, we   record the following special cases of Theorem~\ref{thm0}.
\begin{corollary}\label{cor0}
Let $\psi$ be a  nonzero entire function  on  $\CC$ and $0<p,q\leq  \infty$.
\begin{enumerate}
\item If $p\leq q$, then  the  statements:
$ S_{\psi}(\mathcal{F}_p)\subseteq  \mathcal{F}_q$,
 $\psi(z)= az+b$ for some $a, b \in \CC$,
 $ S_{\psi}:  \mathcal{F}_p\rightarrow  \mathcal{F}_q$   is bounded and  $ S_{\psi}:  \mathcal{F}_p\rightarrow  \mathcal{F}_q$ is globally Lipschitz continuous,  are all equivalent.
\item If $p>q$, then the  statements: $ S_{\psi}(\mathcal{F}_p)\subseteq  \mathcal{F}_q$,  $\psi=$ constant,  $ S_{\psi}:  \mathcal{F}_p\rightarrow  \mathcal{F}_q$ is bounded, and
$ S_{\psi}:  \mathcal{F}_p\rightarrow  \mathcal{F}_q$ is globally Lipschitz continuous,  are all equivalent.
\item $ S_{\psi}:  \mathcal{F}_p\rightarrow  \mathcal{F}_q$ cannot be compact for  any pair of  $p$ and $q$.
\end{enumerate}
\end{corollary}
 Given the fact that  $\mathcal{F}_p$ properly contains the space $\mathcal{F}_q$ when $q<p$ \cite[Theorem 2.10]{Zhu}, it should be clear that    a superposition from the former to the latter is possible only
via constant functions.
 \begin{corollary}
 Let $u$  be a   nonzero   entire function on $\CC$ and $0<p, q\leq \infty$.  Then
 \begin{enumerate}
 \item if $p\leq q$, then   the statements: $ M_u ( \mathcal{F}_p)\subseteq \mathcal{F}_q$,  $u= $   constant,  $ M_u:  \mathcal{F}_p\rightarrow  \mathcal{F}_q$ is bounded,  and
 $ M_u:  \mathcal{F}_p\rightarrow  \mathcal{F}_q$ is globally Lipschitz continuous, are all equivalent.
 \item if $p>q$, then
 $ M_u$  fails to map  $\mathcal{F}_p$ into $  \mathcal{F}_q$.
 \item $ M_u:  \mathcal{F}_p\rightarrow  \mathcal{F}_q$ cannot be compact for  any pair of $p$ and $q$.
      \end{enumerate}
 \end{corollary}
We close  this section with a word on notation.  The notion
 $U(z)\lesssim V(z)$ (or
equivalently $V(z)\gtrsim U(z)$) means that there is a constant
$C$ such that $U(z)\leq CV(z)$ holds for all $z$ in the set of a
question. We write $U(z)\simeq V(z)$ if both $U(z)\lesssim V(z)$
and $V(z)\lesssim U(z)$.
\section{Proof of  the  results}
Let us first consider the following  necessity condition  about  the superposition operator $ S_{\psi}$. The proof of the lemma  demonstrates  an important test function
that will be used in the proof of the main result.
\begin{lemma}\label{lem1}
Let $\psi$ be  nonzero  entire function on $\CC$ and $0<p,  q \leq \infty$. If $ S_{\psi}$ maps $ \mathcal{F}_p$  into $  \mathcal{F}_q$, then  $\psi= az+b$ for some $a, b\in \CC$.
\end{lemma}
\begin{proof}
 To prove the assertion, we need to choose test function from $ \mathcal{F}_p$ for any fixed  $p$.  Thus, for  $c$ with $\frac{1}{4}<c<\frac{1}{2}$, consider the function
\begin{align*}f_c(z)= e^{cz^2}\end{align*}
In view of the estimate,
\begin{align*}
\big| f(z)e^{-|z|^2/2}\big|\leq  e^{(c-\frac{1}{2})|z|^2},
\end{align*}
 the functions $f_c\in \mathcal{F}_p$. Moreover, the Fock space norm is invariant under rotations so if $|\lambda| = |\mu| = 1$,  then all functions
\begin{align*}
f_{\lambda, \mu,c} (z)= \lambda e^{c\mu z^2} \in \mathcal{F}_p
\end{align*}and have the same norm as $f_c$ which can be seen by making change of variables $\tau= \sqrt\mu z$ in the integration  for $p\neq \infty$ and in  estimating  the supremum norm  for $p= \infty$. By the same reasoning, we also observe  that
\begin{align*}
\|S_\psi f_{\lambda, \mu,c}\|_q= \|S_\psi f_{c}\|_q.
\end{align*}
Now for any $w\neq 0$, we can find   $z, \lambda, \mu \in \CC$  with $|\lambda| = |\mu| = 1$ and
\begin{align*}
w=\lambda e^{c\mu z^2}
\end{align*} and also that $c\mu z^2$ is real and positive  which implies $w/\lambda= \overline{\lambda}w $. Clearly, for large $w$, we have
\begin{align*}
|z|^2
= \frac{1}{c} \log( \overline{\lambda}w)\simeq \frac{\log |w|}{c}.
\end{align*}
Therefore, using   the estimate above and  in \eqref{global}
\begin{align*}
|\psi(w)|= |\psi\big(\lambda e^{c\mu z^2}\big)| \leq e^{|z|^2/2}\|S_\psi f_{\lambda, \mu,c}\|_q= e^{|z|^2/2} \|S_\psi f_{c}\|_q= \|S_\psi f_{c}\|_q |w|^{1/2c}
\end{align*} for sufficiently large $w$.\\
 Since  $\|S_\psi f_{c}\|_q$ is  a constant value and $1/2c <2$, the standard Cauchy
estimates imply that $\psi$ is a polynomial of degree at most one as desired.
\end{proof}
Next,  we recall the notion of order and type of an entire function and prove one more important lemma.
Let $f$ be an entire function, and $ M(r, f) =\max_{|z|=r}|f(z)|$.  The order of $f$ is
\begin{align*}
\rho(f) = \limsup_{r\to \infty} \frac{\log \log M(r,f)}{\log r}.
\end{align*}
 If $0<\rho<\infty$, then the type of  $f$ is given by
 \begin{align*}
 \tau(f)= \limsup_{r\to \infty} \frac{\log M(r,f)}{ r^{\rho(f)}}.
 \end{align*}
By definition, if a function $f$ has order $\rho$, then for every $\epsilon >0$
\begin{align*}
M(f,r)= O(e^{\rho+\epsilon})
\end{align*} when $r\to \infty$. This clearly shows that any function of order less than 2 belongs to all the Fock spaces $\mathcal{F}_p$.  Conversely, by  \cite[Theorem 2.12]{Zhu},  every function in
$\mathcal{F}_p$  has order at most 2, and if its order is exactly  2, it must be of type less than or equal to $1/2$. In the next lemma, we prove that if the function has no zeros, then its type can not be $1/2$.
\begin{lemma}\label{order} Let $0<p\leq \infty$ and $u$ is a non-vanishing function in  $\mathcal{F}_p$. Then \begin{align}
 u(z)= u(0) e^{a_1z+a_2z^2}, \ \  a_1, a_2 \in \CC \ \text{ and}\   |a_2|< 1/2.
 \end{align}
\end{lemma}
\begin{proof} We argue as follows. Since $u\in \mathcal{F}_p$, by \cite[Theorem 2.12]{Zhu}, $ \rho(u)\leq 2$  and if $ \rho(u)= 2$, it must be of type less than or equal to $1/2$. On the other hand, since $u$ is non-vanishing, it follows from the Hadamard Factorisation Theorem that
 \begin{align}
 u(z)=  e^{a_0+a_1z+a_2z^2}= u(0) e^{a_1z+a_2z^2}
 \end{align} where $ a_0, a_1, a_2 \in \CC$ and $|a_2|\leq 1/2$. It remains to show that $|a_2|=1/2 $ can not happen. Assuming to the contrary, we may simply set $a_2= 1/2$ (if not since the Fock space norm is invariant under rotation, we can find a $\mu$ with  $|\mu|=1$  such that $\mu a_2= 1/2 $ ).
 Then for $p<\infty$
 \begin{align*}
 \int_{\CC} |u(z)|^p e^{-\frac{p}{2}|z|^2} dA(z)=  |u(0)|^p \int_{\CC} e^{p\big(\Re(a_1z)+1/2\Re(z^2))-\frac{p}{2}|z|^2} dA(z)\\
 =|u(0)|^p  \int_{\R}\int_{\R}e^{p\Re(a_1)x-p(\Im(a_1)y+ y^2)} dx dy\\
 = |u(0)|^p \bigg( \int_{\R}e^{p\Re(a_1)x} dx\bigg) \bigg( \int_{\R} e^{-p(\Im(a_1)y+ y^2)}dy\bigg)= \infty
 \end{align*} as the first integral with respect to $x$ diverges while the second integral with respect to $y$ converges which contradicts the assumption that $u\in \mathcal{F}_p$.\\
 If $p= \infty$, we replace the above integral argument with supremum to arrive at the same conclusion.
 \end{proof}
 For the sake of better exposition, we split Theorem~\ref{thm0} into two theorems below.
   \begin{theorem}\label{thm3}
Let $\psi$ and $u$ be nonzero entire functions on  $\CC$ and $0<p<q\leq  \infty$. Then if
\begin{enumerate}
 \item $p\leq q$, then $ S_{(u,\psi)}$ maps $\mathcal{F}_p$ into $  \mathcal{F}_q$ if and only if
either  $\psi(z)= az+b$ for some $a, b \in \CC$  and $u$ is a constant or $\psi$ is a constant and $u \in  \mathcal{F}_q$. In  the case when $u $ is non-vanishing, it has the form
\begin{align}
\label{form1}
 u(z)= u(0) e^{a_1z+a_2z^2}, \ \  a_1, a_2 \in \CC \ \text{ and}\   |a_2|< 1/2.
 \end{align}
 \item  $p>q$, then $ S_{(u,\psi)}$ maps $\mathcal{F}_p$ into $  \mathcal{F}_q$ if and only if
 $\psi$ is a constant and $u\in \mathcal{F}_q$. In  the case when $u $ is non-vanishing, it has the form in \eqref{form1}.
 \end{enumerate}
\end{theorem}
 \subsection*{Proof of part (i) }
The sufficiency of the condition is  easy to verify. It is the necessity  that may require a  new techniques which we present below.
We shall first prove that if $S_{(u,\psi)}$ maps  $\mathcal{F}_p$ into $ \mathcal{F}_q$, then $ \psi(z)= az+b$ for some $ a, b\in \CC$. The assumption implies $u\in  \mathcal{F}_q$ as already seen before. To this end, if  $u$ is non-vanishing, by   Lemma~\ref{order}
\begin{align}
\label{form102}
 u(z)= u(0)e^{a_1z+a_2z^2}, \  a_1, a_2 \in \CC \ \text{ and}\   |a_2|< 1/2.
 \end{align}
    Suppose for the purpose of contradiction that $ \psi$ is not linear. Then there exists a sequence  $w_n\in \CC$ such that $|w_n| \to \infty$ as $n\to \infty$ and $|\psi(w_n) |\geq  n|w_n|^2$ for  $n\in\NN$. Arguing as in the proof of Lemma~\ref{lem1}, for each $w_n\neq 0$, we can find   $z_n, \lambda_n, \mu_n \in \CC$  with $|\lambda_n| = |\mu_n| = 1$ and
\begin{align}
\label{zero}
w_n=\lambda_n e^{c\mu_n z_n^2}
\end{align} and also that $c\mu_n z_n^2$ is real and positive.  We may further   pick a sparse  subsequence of   $z_n$  such that the discs $D(z_n, 1)$ are mutually disjoint. Now, for  $q<\infty$, applying the operator to the sequence \begin{align*}
f_{c,\lambda_n, \mu_n}(z)= \lambda_n e^{c\mu_n z^2}\end{align*}
 and eventually  invoking the estimate in \eqref{global},
\begin{align}
\label{finite}
 \int_{\CC} |S_{u,\psi} f_{c,\lambda_n, \mu_n}(z)|^q e^{-\frac{q}{2}|z|^2} dA(z)= \int_{\CC} |u(z)|^q |\psi(\lambda_n e^{c\mu_n z^2}|^q e^{-\frac{q}{2}|z|^2} dA(z) \quad  \nonumber\\
 \geq \sum_{n=1}^\infty  \int_{D(z_n,1)} |u(z)|^q|\psi(w_n)|^q e^{-\frac{q}{2}|z|^2} dA(z)\nonumber\\
 \gtrsim \sum_{n=1}^\infty  |u(z_n)|^q|\psi(w_n)|^q e^{-\frac{q}{2}|z_n|^2}
 \geq\sum_{n=1}^\infty n^q |u(z_n)|^q  e^{2cq\mu_n z_n^2-\frac{q}{2}|z_n|^2}.
 \end{align}
  We consider two separate cases depending on whether the multiplier function $u$ has zeros or not in $\CC$.  If $u$  is non-vanishing, by \eqref{form102} and \eqref{finite}
 \begin{align}
\label{finite2}
 \int_{\CC} |S_{u,\psi} f_{c,\lambda_n, \mu_n}(z)|^q  n^q e^{-\frac{q}{2}|z|^2} dA(z) \gtrsim \sum_{n=1}^\infty  |e^{a_0+a_1z_n+a_2z_n^2}|^q |e^{2cq\mu_n z_n^2-\frac{q}{2}|z_n|^2}.
 \end{align}
 Now, if $a_2=0$, then we may chose $c$  in the interval   $\big(\frac{3}{8},\frac{1}{2})$ such that $2c-\frac{1}{2}>0$.  Thus, the last sum in \eqref{finite} diverges. On the other hand, if  $0<|a_2|,$ then  we may chose $c$ in  $\big(\frac{1}{4} + \frac{|a_2|}{2},\frac{1}{2})$ such that
 \begin{align*} 2c-\frac{1}{2} > 2\bigg(\frac{1}{4} + \frac{|a_2|}{2}\bigg)-\frac{1}{2}= |a_2|>0\end{align*}  and hence the sum in \eqref{finite} still diverges. This is a contradiction and hence  $\psi(z)= az+b$ for some $a, b\in \CC$ in this case.\\
Next, assume  that $u$ has zeros and analyze the case when the sequence $z_n$ in the sum \eqref{finite} could belong to the zero set of $u$. Let us for example assume that $u$ is a polynomial. Then there exist positive constants $C$ and $R$ such that for $|u(z)|\geq C$  for $|z|\geq R$. It follows that
\begin{align*}
\|S_{(u,\psi)}f\|_q^q=
 \frac{q}{2\pi}\int_{\CC} \big|u(z) S_\psi f(z)\big |^q e^{-\frac{q}{2}|z|^2} dA(z)\quad \quad \quad \quad \quad \quad \quad \quad\\
  \geq C^q\frac{q}{2\pi}\int_{\{z\in \CC: |z|\geq R\}}\big| S_\psi f(z)\big |^q  e^{-\frac{q}{2}|z|^2} dA(z).
\end{align*} By Lemma~\ref{lem1}, the last  integral above is finite only when $\psi(z)= az+b$.\\
Now  if infinitely many $z_n's$ belong to the zero set of $u$ such that last sum in \eqref{finite} converges, then we keep changing our choice of the  constant  $c$ in the interval  $(1/4,1/2)$ above until  the convergence is not possible. We explain how this can done below.  From \eqref{zero} and  since $c$ has uncountably  many possible values in the said interval, we collect all these values to see that
\begin{align*}
\bigg\{\frac{1}{\sqrt{c}}\frac{\log (w_n/\lambda_n)}{\mu_n}: (c,n) \in (1/4,1/2)\times \NN\bigg\}
\end{align*} is uncountable.  As known, an entire function can not have uncountable zero set. Taking this into account, we  can   choice   the sequence $z_n$ in such a way that it does not belong to the zero set of $u$ when $n\to \infty$.\\
  Now, if  $u$ has order less than $2$, then every choice of  $c$ in the interval $(1/4,1/2)$ gives that the sum in \eqref{finite} diverges. On the other hand,  if $\rho(u)=2$ and  its type, $\tau(u)$,  is less than $1/2$, then we can still choice
$c$ in the interval $( \frac{\tau(u)}{2}+\frac{1}{4}, \frac{1}{2})$ such that sum in \eqref{finite}  diverges.
Thus,  it remains to consider  the extremal case  when
\begin{align}
\label{setting}
|u(z_n)| \simeq e^{-\frac{1}{2}|z_n|^2}
\end{align} as $n\to \infty$. To this end,  setting the expression in \eqref{setting} in \eqref{finite}
\begin{align}
\label{finitee}
 \int_{\CC} |S_{u,\psi} f_{c,\lambda_n, \mu_n}(z)|^q e^{-\frac{q}{2}|z|^2} dA(z)
 \gtrsim \sum_{n=1}^\infty n^q |u(z_n)|^q  e^{2cq\mu_n z_n^2-\frac{q}{2}|z_n|^2}\nonumber\\
 \simeq \sum_{n=1}^\infty n^q |u(z_n)|^q  e^{2cq\mu_n z_n^2-q|z_n|^2}
 \end{align} holds for   all  $c$  in  $(1/4,1/2)$.  Letting $c\to \frac{1}{2}$, we note that the sum in \eqref{finitee} converges only when  the sequence $\{n^q:\in \NN\}$
 is summable,  which is a contradiction again, and   hence $\psi(z)= az+b$ for some constants $a, b\in \CC $.\\
 Next, we   set $\psi(z)= az+b$ and  show  that $u$ is necessarily  a constant function whenever   $a\neq 0$. Aiming to argue in the contrary, assume that $u$ is not a constant. Then we can pick a sparse sequence   $w_n$  such that $|u(w_n)| \to \infty$ as $n\to \infty$ and that the discs $D(w_n, 1)$ are mutually disjoint. For each positive $\epsilon$, observe that the  function  $g_\epsilon(z)= e^{\frac{z^2}{2+\epsilon}}\in  \mathcal{F}_p$ for all $p$.  Then we can find a sequence  $\mu_n$ such that $|\mu_n|=1$ and $\mu_n w_n^2 $ are real  and positive for sufficiently large $n$.      Setting  $ h_{\epsilon,\mu_n} (z)= e^{\frac{\mu_nz^2}{2+\epsilon}}-\frac{b}{a}$ and  eventually  applying  \eqref{global}
\begin{align*}
\|S_{(u,\psi)} h_{\epsilon,\mu_n}\|_q^q=
 \frac{q}{2\pi}\int_{\CC} \big|S_{(u,\psi)}h_{\epsilon,\mu_n}(z)\big |^q e^{-\frac{q}{2}|z|^2} dA(z)\quad \quad\quad \quad\quad \quad\quad \quad\\
 = \frac{q}{2\pi}|a|^q\int_{\CC} \big|u(z)  h_{\epsilon,\mu_n}(z)\big |^q e^{-\frac{q}{2}|z|^2} dA(z)\\
\geq \frac{q}{2\pi} |a|^q\sum_{n=1}^\infty \int_{D(w_n, 1)} \big|u(z)  h_{\epsilon,\mu_n}(z)\big |^q e^{-\frac{q}{2}|z|^2} dA(z)  \ \ \ \quad\\
\geq \frac{q}{2\pi} |a|^q\sum_{n=1}^\infty |u(w_n)|^{q}  e^{\big(\frac{q}{2+\epsilon}-\frac{q}{2}\big)w_n^2} \ \ \ \quad \quad
\end{align*} for all   $\epsilon>0$.  Letting $\epsilon \to 0$, we observe that the sum diverges,  which  is again  a contradiction.\\
For  $q= \infty,$  we may simply replace the integral above by the supremum norm and argue to arrive at the same conclusion. This completes the proof of part (i).
\subsection*{Proof of part (ii)}
The sufficiency of the condition is  clear again.
Let $p>q$. Then as  already proved above in the first part, $\psi(z)= az+b$ is a necessary condition for $ S_{(u,\psi)}$  to map  $ \mathcal{F}_p$ into $ \mathcal{F}_q$ independent of the size of $p$ and $q$. Assume that $a\neq 0$ and for any $f\in \mathcal{F}_p $, consider the function $f_b= f-\frac{b}{a}\in \mathcal{F}_p$ and observe
 \begin{align*}
 \big\| S_{(u,\psi)} f_b\big \|_q = |a| |\alpha|\|f\|_q
   \end{align*} Then,  our conclusion follows from the fact that $\mathcal{F}_p\backslash \mathcal{F}_q $ is non-empty \cite[Theorem 2.10]{Zhu}   and $|a| |\alpha|\neq 0$ and  completes the proof of  Theorem~\ref{thm3}.
 \begin{theorem}\label{thm4}
Let $\psi$ and $u$ be  nonzero   entire functions on $\CC$ and $0<p,q\leq \infty$. If
$ S_{(u,\psi)}$ maps $\mathcal{F}_p$ into $ \mathcal{F}_q$, then it is bounded and globally   Lipschitz continuous.
But  $S_{(u,\psi)}$ cannot be compact.
 \end{theorem}
\emph{Proof}.
    Clearly, f $\psi$ is a constant and $u\in \mathcal{F}_q $, then  $ S_{(u,\psi)}: \mathcal{F}_p\to   \mathcal{F}_q$ is bounded. Thus, lets dispose the sufficiency of the condition when   $\psi(z)= az+b$ and $u= \alpha$. Consider   $ f\in \mathcal{F}_p$  and compute
\begin{align*}
\|S_{(u,\psi)} f\|_q^q=
 \frac{q}{2\pi}\int_{\CC} \big|u(z)(a f(z)+b)\big |^q e^{-\frac{q}{2}|z|^2} dA(z) \ \ \ \quad \quad  \quad \quad  \quad \quad \quad \quad\\
\leq \frac{q}{2\pi} |2\alpha|^{q}\Big(|a|^q \int_{\CC} | f(z) |^q e^{-\frac{q}{2}|z|^2} dA(z)+ \int_{\CC}|b|^qe^{-\frac{q}{2}|z|^2}dA(z)\Big)\\
= |2\alpha|^{q}|a|^q \|f\|_q^{q}+ |2\alpha|^{q}|b|^q \leq  |2\alpha|^{q}|a|^q \|f\|_p^{q}+ |2\alpha|^{q}|b|^q <\infty
\end{align*} where the second  inequality follows from \cite[Theorem 2.10]{Zhu}.

Next, let us show that the spaces support no compact weighted  superposition operators. Aiming to arrive at a contradiction, let  $ S_{(u,\psi)}: \mathcal{F}_p\to \mathcal{F}_q$  be compact and hence $\psi(z)= az+b$ for some  complex number $a, b$ and $u= \alpha$ or  $\psi= b$ and $u\in \mathcal{F}_q$. The sequence  $k_w \in \mathcal{F}_p$ is bounded and  converges to zero uniformly on compact subsets of $\CC$ when $|w|\to \infty$.  Applying the operator to  $k_w $ and eventually  \eqref{global}
 \begin{align*}
 \big\| S_{(u,\psi)} k_w\big \|_q \geq  |ak_w(w)+b||\alpha| e^{-|w|^2/2}=  |a+b e^{-|w|^2/2}||\alpha|
 \end{align*}
  and $ \big\| S_{(u,\psi)} k_w\big \|_q   \to 0$ as $ \ |w| \to \infty
 $ only when $a= 0$.  It follows that
 \begin{align*}
 \big\| S_{(u,\psi)} k_w\big \|_q =   |b||\alpha| \|1\|_q \to 0
 \end{align*} as $ \ |w| \to \infty $ only if  $b=0$ which contradicts that $\psi$ is nonzero.\\
   If $a=0$ and $u\in \mathcal{F}_q$ is non-zero, the same conclusion follows easily.

It remains to show that $S_{(u,\psi)}: \mathcal{F}_p\to \mathcal{F}_q$ is globally  Lipschitz continuity.
   Since the case for $p>q$ is trivial, we assume  $p\leq q$. An application of Theorem~\ref{thm3} gives  that  $\psi(z)= az+b$ and $u= \alpha\neq0$ or $\psi= b$ and $u\in \mathcal{F}_q$.  If $a=0$, then the conclusion follows easily. Thus,  assume $a\neq0$ and hence $u = \alpha\neq0$.    For   $f, g \in \mathcal{F}_p$
\begin{align*}
\|S_{(u,\psi)} f-S_{(u,\psi)} g\|_q^q= \|auf- aug\|_q= |a\alpha|\|f- g\|_q\leq  |a\alpha |\| f-g\|_p
 \end{align*} where the last inequality follows from the inclusion property.  This shows that  $ S_{(u,\psi)}$ is  globally   Lipschitz continuous with constant $|a\alpha|$.
\subsection{Remark} In the remainder of this section, we present  alternative proofs for the necessity  in Theorem~\ref{thm3} when
 \begin{enumerate}
\item \   $q= \infty$
  \item \ $2=p<q<\infty$
  \item \   $q<p<\infty$ and $u$ is non-vanishing.
  \end{enumerate}   The proof is interest of its own as it shows  how  the zero sets and  uniqueness sets
in Fock spaces play important roles in the study of weighted superposition operators. Unfortunately, a complete characterization of such
sets are still an open problem as far as we know. Some necessary and sufficient conditions can be read in \cite[Chapter 5]{Zhu}.\\
 (i)  We consider the square lattice in the complex plane
 \begin{align*}
 \Lambda= \big\{\omega_{mn}=\sqrt{\pi} (m+in): m, n\in \Z \big\}
 \end{align*} where $\Z$ denotes  the set of all integers. As known the Weierstrass $\sigma$-function associated to  $ \Lambda$ is defined by
 \begin{align*}
\sigma(z)= z\prod_{\substack{(m, n)\neq (0,0)\\ (m, n) \in \Z^2}} \bigg(1-\frac{z}{\omega_{mn}}\bigg)\exp\bigg( \frac{z}{\omega_{mn}}+ \frac{z^2}{2\omega_{mn}^2}\bigg).
 \end{align*}
 Furthermore, it is known that $ \Lambda$ is the zero set for  $\sigma$ and by \cite[Lemma 5.6]{Zhu}, $\sigma \in \mathcal{F}_\infty $ but not in  any of the other Fock spaces $\mathcal{F}_p$.
 Aiming to argue in the contrary, suppose now that $u\neq 0$ and $\psi$ is not a constant.  Since $ S_{(u,\psi)}$ maps $\mathcal{F}_\infty $ into $  \mathcal{F}_q$,  the function
 \begin{align*}
 F= S_{(u,\psi)} \sigma-\psi(0) u= u \psi(\sigma)-\psi(0)u \in \mathcal{F}_q.
 \end{align*}
 Now since $\psi$ is not a constant, the function $F$ is non-zero and vanishes on $ \Lambda$. On the other hand, by \cite[Lemma 5.7]{Zhu},   $ \Lambda$ is a uniqueness set for  $\mathcal{F}_q$  for all $q\neq \infty$ which implies that
 $F=0$. This is contradicts the fact that $F\neq0$. Therefore, $\psi$ is a constant.

 (ii)  We consider the set $\Omega= \Lambda-\{0\}$. As shown in \cite[p.204]{Zhu}, $\Omega$ is a uniqueness set for  $\mathcal{F}_2$.  Indeed,  the function $g(z)= \sigma(z)/z$ belongs to $\mathcal{F}_p$ if and only if
 $p>2$. Now to argue as above, assume $\psi$ is not a constant and  hence
 \begin{align*}
 G(z)=  S_{(u,\psi)} g-\psi(0) u= u \psi(g)-\psi(0)u \in \mathcal{F}_2
 \end{align*} is non-zero and vanishes on $\Omega$, but $\Omega$ is a uniqueness set for $\mathcal{F}_2$. Thus,  $G=0$ resulting a  contradiction again.

 (iii) For a positive number $R$, we consider the  following modified lattice
 \begin{align*}
 \Lambda_{R}= \big\{\omega_{mn}: m, n\in \Z \big\}
 \end{align*} where
 \begin{align*}
 \omega_{mn}= \begin{cases}
 & z_{mn} ,\  \text{if} \  n\neq0 \ \text{or}\ n= m=0\\
 & \sqrt{\pi} \big( m+\frac{Rm}{|m|}\big),  \  \text{if} \  n=0 \ \text{and }\  m\neq0
 \end{cases}
 \end{align*}
 Then  the modified Weierstrass function associated to $\Lambda_{R}$ is given by
  \begin{align*}
\sigma_{R}(z)= z\prod_{\substack{(m, n)\neq (0,0)\\ (m, n) \in \Z^2}} \bigg(1-\frac{z}{\omega_{mn}}\bigg)\exp\bigg( \frac{z}{\omega_{mn}}+ \frac{z^2}{2z_{mn}^2}\bigg).
 \end{align*}
Now if we choose  $R$  to be a number such that $
\frac{1}{p} <R< \frac{1}{q},
$ then as shown in the proof of \cite[Theorem 1.1]{ABO}, the function $\sigma_{R}\in \mathcal{F}_p$ and $\sigma_{R}$ constitutes a zero set  for  $\mathcal{F}_p$ while it fails to be
a zero set for $\mathcal{F}_q$.\\
Now to argue as in the previous two cases, assume $u\in \mathcal{F}_q$ is non-vanishing and $\psi$ not a constant.  Then
\begin{align*}
 H= S_{(u,\psi)} \sigma_{a,R}-\psi(0) u= u \psi(\sigma)-\psi(0)u \in \mathcal{F}_q
 \end{align*} is non constant  and  vanishes only on the set $\Lambda_{R}$. It follows that $\Lambda_{R}$ is a zero set for $\mathcal{F}_q$ which is a contradiction.

\end{document}